\documentclass{amsart}
\usepackage{amssymb,amsfonts,latexsym}

\usepackage{times}
\numberwithin{equation}{section}

\newtheorem{theorem}{Theorem}[section]
\newtheorem{lemma}[theorem]{Lemma}
\newtheorem{proposition}[theorem]{Proposition}
\newtheorem{corollary}[theorem]{Corollary}
\theoremstyle{definition}
\newtheorem{definition}[theorem]{Definition}

\newtheorem{example}[theorem]{Example}

\newtheorem{remark}[theorem]{Remark}

\newcommand{\Rep}{\text{Rep}}

\newcommand{\g}{\mathfrak{g}}
\newcommand{\h}{\mathfrak{h}}

\newcommand{\ot}{\otimes}

\newcommand{\ben}{\begin{enumerate}}
\newcommand{\een}{\end{enumerate}}

\newcommand{\Vect}{\text{Vec}}

\newcommand{\CC}{{\mathbb{C}}}

\newcommand{\cC}{{\mathcal C}}

\begin{document}

\title[Quasisymmetric and unipotent tensor categories]
{Quasisymmetric and unipotent tensor categories}
\author{Pavel Etingof}
\address{Department of Mathematics\\ Massachusetts Institute of
Technology\\
Cambridge, MA 02139, USA} \email{etingof@math.mit.edu}

\author{Shlomo Gelaki}
\address{Department of Mathematics\\ Technion-Israel Institute of
Technology\\ Haifa 32000, Israel}
\email{gelaki@math.technion.ac.il}

\date{October 31, 2007}

\maketitle

\section{Introduction}

One of the most important early developments in the theory of
quantum groups was Drinfeld's classification, in characteristic zero,
of quasitriangular quasi-Hopf QUE (quantized universal enveloping)
algebras \cite{Dr1,Dr2}. In the language of tensor categories,
this is, in essence, a classification of 1-parameter flat formal deformations,
as a braided category\footnote{Here by a
flat formal deformation of the category of
finite dimensional modules over a (topological) algebra $A_0$ we
mean the category of finite dimensional modules over a flat formal
deformation $A$ of $A_0$. This definition will suffice for our
purposes; we note, however, that in general,
the definition of a flat formal deformation of an abelian
category is fairly nontrivial (\cite{LV}).},
of the representation category of a Lie
algebra\footnote{Throughout the paper, we work
over the ground field $\mathbb C$ of complex numbers.} $\g_0$.
The answer is that such
deformations are parameterized by pairs $(\g,t)$, where
$\g$ is a flat formal deformation of $\g_0$, and $t$ is
an element in $(S^2\g)^\g$. More specifically, the deformed
category is the category of representations of $\g$,
with the usual tensor product functor,
the braiding is given by the formula $\beta=P\circ
e^{\hbar t/2}$, where $P$ is the flip map,
and the associativity isomorphism is $\Phi(\hbar t_{12},\hbar
t_{23})$, where $\Phi(a,b)$ is any Drinfeld associator.
In particular, {\it symmetric} deformations correspond
to the case $t=0$; in other words, such deformations
come simply from deformations $\g$ of the corresponding Lie
algebra $\g_0$.

Drinfeld's result generalizes mutatis mutandis to a more general
setting where $\g_0$ is a Lie superalgebra, and to the situation
when $\g_0$ is replaced with an affine proalgebraic supergroup
$G_0$. More specifically, suppose that $G_0$ is an affine
proalgebraic supergroup, and
$u_0\in G_0$ an element of order 2 acting by parity on the
function algebra ${\mathcal O}(G_0)$. Let
$\Rep(G_0,u_0)$ be the category of representations of $G_0$ on
finite dimensional supervector spaces, in which $u_0$ acts by
parity. Then Drinfeld's work implies that any flat formal
deformation of the category $\Rep(G_0,u_0)$, as a braided
category, has the form $\Rep(G,u)$, where $G$ is a deformation of
$G_0$, $u$ is the (unique) deformation of $u_0$ in $G$, and the
associativity isomorphism and braiding are $\Phi(\hbar
t_{12},\hbar t_{23})$ and $P\circ e^{\hbar t/2}$, for some $t\in
(S^2\g)^G$ (here $\g={\rm Lie}(G)$). Moreover, $(G,u,t)$ are
determined by the deformation uniquely up to an isomorphism. In
particular, if the deformation is symmetric then $t=0$, and the
deformation is $\Rep(G,u)$.

By Deligne's theorem \cite{De2}, any
symmetric tensor category of exponential growth
has the form $\Rep(G_0,u_0)$; thus, Drinfeld's result
provides a description of flat formal deformations of any
 symmetric tensor category with at
most exponential growth. In particular, in the special case
of symmetric deformations, Drinfeld's result can be viewed as
a formal analog of Deligne's theorem.

Unfortunately, Drinfeld's method makes a serious use of the
presence of the formal parameter $\hbar$, i.e., of the
fact that the braided categories at hand are symmetric
modulo this parameter. For this reason, it cannot be applied
to classifying braided categories over $\CC$, not involving
$\hbar$ (even those of exponential
growth). In fact, we are very far from the classification of
such categories, even in the special case of finite semisimple
(i.e., fusion) categories.

On the other hand, it turns out that
there is a subclass of braided categories
for which Drinfeld's method does work. This is the class of
{\it quasisymmetric} categories, introduced essentially in
\cite{EK2}\footnote{To be precise, the definition of a
quasisymmetric category in \cite{EK2} is
somewhat different from the one used in this paper, but the difference
is inessential for what we do.}. They are, by definition,
 braided categories with exponential
growth, in which the square of the braiding is the identity
on the product of any two simple objects. For such categories,
the infinite power series in $\hbar$ which occur in Drinfeld's
construction terminate when applied to the tensor product of
any objects, and thus become polynomials in $\hbar$. These
polynomials can then be evaluated at $\hbar=1$, which
allows one to apply Drinfeld's method to the situation without
$\hbar$.

The study of quasisymmetric categories by means of Drinfeld's
method is one of the main goals of this paper. Another is to study
{\it unipotent} tensor categories, i.e., such that every simple
object is the neutral object. Namely, in Section 2, using
Drinfeld's method and Deligne's theorem, we give a classification
of quasisymmetric categories, which is similar to Drinfeld's
classification of quasitriangular quasi-Hopf QUE algebras.
Specifically, we show that equivalence classes of such categories
are in bijection with equivalence classes of triples $(G,u,t)$,
where $G$ is an affine proalgebraic supergroup, $u\in G$ an
element of order 2 which acts by parity on the function algebra
${\mathcal O}(G)$, and $t$ a nilpotent element of $(S^2\g)^G$,
where $\g={\rm Lie}(G)$. This gives a generalization of Deligne's
theorem to the case of quasisymmetric categories. As a special
case, this result yields a classification of braided unipotent
categories.  In Section 3, we proceed to classify unipotent fiber
functors on quasisymmetric categories, i.e., functors that are
``standard'' on the canonical symmetric part of the category (the
subcategory $\otimes$-generated by the simple objects). Namely, we
show that such functors are in bijection with nilpotent solutions
$r$ of the classical Yang-Baxter equation, such that $r+r^{21}=t$.
In Section 4, using the quantization theory from
\cite{EK1,EK2,EK3}, we classify coconnected Hopf algebras (i.e.,
Hopf algebras with a unique simple comodule) by showing that they
all come from quantization of prounipotent Poisson
proalgebraic groups. This provides a classification of unipotent
tensor categories with a fiber functor.

\section*{Acknowledgements}
The research of the first author was
partially supported by the NSF grant DMS-0504847.
The second author was supported by The Israel
Science Foundation (grant No. 125/05). He also thanks MIT for its
warm hospitality during his Sabbatical.
Both authors were supported by BSF grant No. 2002040.

\section{Quasisymmetric categories}

Let $\mathcal{C}$ be a rigid tensor category over $\mathbb{C}$. In
particular, $\cC$ is an abelian category, with finite dimensional
spaces of morphisms and all objects having finite length, and we
have ${\rm End}({\bold 1})=\mathbb C$. We will also assume
throughout the paper that
 $\cC$ has exponential growth, i.e.,
that for every object $Y$ there exists $d(Y)\ge 1$ such that
${\rm length}(Y^{\otimes n})\le d(Y)^n$ for all $n\ge 1$.

\subsection{Unipotent categories}

\begin{definition}
The tensor category $\mathcal{C}$ is called unipotent if the only simple object
in $\mathcal{C}$ is the neutral object ${\bf 1}$.
\end{definition}

\begin{remark} Note that the exponential growth condition
is automatic for unipotent tensor categories.
\end{remark}

The simplest example is the following one.

\begin{example}
Let $G$ be a prounipotent proalgebraic group. Then
$\mathcal{C}:=\Rep(G)$, the category of rational representations
of $G$, is unipotent.
\end{example}

\begin{proposition}\label{symm}
If $\mathcal{C}$ is a symmetric unipotent category then
$\mathcal{C}=\Rep(G)$ for some prounipotent proalgebraic
group $G$.
\end{proposition}

\begin{proof}
This follows from Deligne's theorem \cite{De1} since the categorical dimension
of an object $X$ is just its length, hence a positive integer.
\end{proof}

\subsection{Unipotent radical of a supergroup}

Let $G$ be an affine proalgebraic supergroup
\footnote{For a short introduction to supergroups and Deligne's
theorem \cite{De2} we refer the reader, for example, to \cite{EG2}.}.
Denote by $U\subset G$ the intersection of the kernels
of irreducible algebraic representations of $G$, and set $G_{\rm red}:=G/U$.

Thus, we have a natural exact sequence of supergroups
$$
1\to U\to G\to G_{\rm red}\to 1.
$$

\begin{definition}
We will call the supergroup $G_{\rm red}$ the {\it reductive quotient} of
$G$, and $U$ the {\it unipotent radical} of $G$.
We will say that a supergroup $G$ is {\it reductive} if
$U=1$ and $G=G_{\rm red}$.
\end{definition}

Note that for a reductive supergroup, it is not always true that
the category of its representations is semisimple: for example,
the supergroup $GL(m|n)$ is reductive, but it has representations
which are reducible but indecomposable.

\subsection{Quasisymmetric categories}

Let $\cC$ be a braided category with braiding $\beta$.

\begin{definition}
Let us say that $\cC$ is quasisymmetric if for every simple
objects $X,Y\in \cC$, one has $\beta^2={\rm Id}$ on $X\otimes Y$.
\end{definition}

\begin{example}
Every symmetric and every unipotent braided tensor category is
quasisymmetric.
\end{example}

Let $\cC$ be a quasisymmetric tensor category.
Denote by $\cC_s$ the full tensor subcategory of $\cC$
$\otimes$-generated by simple objects of $\cC$ (i.e., formed by
the subquotients of direct sums of tensor products of simple objects).

The following proposition is obvious from the braiding axioms.

\begin{proposition}\label{bra}
The category $\cC_s$ is symmetric.
\end{proposition}

\begin{definition}
We will call $\cC_s$ {\it the canonical symmetric part} of
$\cC$.
\end{definition}

\begin{remark} Note that a quasisymmetric category $\cC$ is unipotent if and
only if $\cC_s$ is the category of finite dimensional vector spaces.
\end{remark}

\begin{example}
Let $G$ be an affine proalgebraic supergroup, and $u\in G$ an element
of order 2 acting by parity on the algebra of regular functions
${\mathcal O}(G)$, and let $\Rep(G,u)$
be the category of representations of $G$ on finite dimensional
supervector spaces on which $u$ acts by parity. Then $\cC_s={\rm
Rep}(G_{\rm red},u)$.
\end{example}

By Deligne's theorem \cite{De2}, Proposition \ref{bra} implies the following corollary.

\begin{corollary}\label{redgru} Let $\cC$ be a quasisymmetric category.
There exists a unique pair $(G_{\rm red},u)$, where $G_{\rm red}$ is a reductive
proalgebraic supergroup, and $u\in G_{\rm red}$ an element of order
2 acting by parity on ${\mathcal O}(G_{\rm red})$, such that
$\cC_s=\Rep(G_{\rm red},u)$.
\end{corollary}

\subsection{Construction of quasisymmetric categories}

Recall that for an affine proalgebraic supergroup $G$,
its Lie (super)algebra ${\rm Lie}(G)$ is defined as the set of left invariant derivations
of ${\mathcal O}(G)$.

Let $G$ be an affine proalgebraic supergroup with Lie algebra
$\mathfrak{g}$. Let $G_{\rm red}$ be the reductive quotient of $G$,
$U$ the unipotent radical of $G$, and
$\mathfrak{g}_r$, ${\mathfrak{u}}$ their Lie algebras.

\begin{definition}
Let us say that an element of the tensor square\footnote{Since $G$
is a affine proalgebraic group, $\g$ is a provector space. Thus
tensor powers of $\g$ will be understood in the completed sense.}
$\g^{\otimes 2}$ is {\it nilpotent} if it projects to zero in
$\g_r^{\otimes 2}$, i.e., if $t\in \g\otimes
{\mathfrak{u}}+{\mathfrak{u}}\otimes\g$.
\end{definition}

Let $u\in G$ be an element of order 2 acting on ${\mathcal O}(G)$
by parity. Let $t\in (S^2\mathfrak{g})^G$
be an invariant nilpotent symmetric $2-$tensor.

Let $\Phi=\Phi(a,b)$ be a Lie associator of Drinfeld;
it is an element in the completed free associative algebra in two
non-commuting variables $a,b$ satisfying some equations \cite{Dr2}
(see also \cite{ES}, page 158).

Now let $\mathcal{C}(G,u,t,\Phi)$ be the braided tensor category defined as follows.
As an abelian category it is just $\Rep(G,u)$. The tensor
product bifunctor is the usual one,
while the associativity constraint is given by
$\alpha:=\Phi(t_{12},t_{23})$ (i.e., $\alpha_{|X\ot Y\ot Z}=
\Phi(t_{12},t_{23})_{|X\ot Y\ot Z}$). The braiding is given by
$\beta:=P\circ e^{t/2}$, where $P$ is the standard flip map.

\begin{remark}\label{nilpo}
Let us explain why $\alpha$ is well defined
(the explanation for $\beta$ is similar).
Recall that any object $X$ in a finite length abelian
category has a canonical filtration $F^\bullet$:
$F^0(X)$ is the sum of all simple subobjects of $X$,
and $F^i(X)$ is defined inductively as the preimage in $X$ of
$F^{i-1}(X/F^0(X))$. It is clear that if $\cC=\Rep(G,u)$, and
$a\in {\mathfrak{u}}$, then $aF^i(X)\subset
F^{i-1}(X)$, so $a$ lowers the filtration degree by $1$.
This implies that if the lengths of $X,Y,Z$ are
$l_X,l_Y,l_Z$, and $n>l_X+l_Y+l_Z-3$, then any product
$t_{i_1j_1}...t_{i_nj_n}$ acts by zero in $X\otimes Y\otimes Z$.
This means that the series defining $\alpha$ terminates, and thus
$\alpha$ is well defined.
\end{remark}

\begin{proposition}\label{phi}
For all $\Phi$, $\Phi'$, $\mathcal{C}(G,u,t,\Phi)$
is equivalent to $\mathcal{C}(G,u,t,\Phi')$.
\end{proposition}

\begin{proof}
Set $\Phi_1:=\Phi(\hbar t_{12},\hbar t_{23})$,
$\Phi_2:=\Phi'(\hbar t_{12},\hbar t_{23})$. By Theorem 3.15 in \cite{Dr1},
there exists an invariant symmetric twist
$T:=T(\hbar t)\in (U(\mathfrak{g})^{\otimes 2}
)^G[[\hbar]]$, given by
a universal formula, such that
$\Phi_1^T=\Phi_2$. Now, since $t$ is nilpotent, similarly to
Remark \ref{nilpo},
$T(\hbar t)_{|X\ot Y}$ is a polynomial in $\hbar$
for all $X,Y$, and hence
can be evaluated at $\hbar=1$. So $T(t)$ is a well defined functorial morphism
$X\ot Y\to X\ot Y$. The identity functor
$\mathcal{C}(G,u,t,\Phi)\to \mathcal{C}(G,u,t,\Phi')$
equipped with the tensor structure $T(t)$ is an equivalence of
braided tensor categories, as desired.
\end{proof}

\begin{remark}
By Proposition \ref{phi}, we may (and will) denote $\mathcal{C}(G,u,t,\Phi)$ simply by
$\mathcal{C}(G,u,t)$. Note that $\cC(G,u,0)=\Rep(G,u)$ as a
braided tensor category. If $u=1$, we will write $\cC(G,t)$ for
$\cC(G,1,t)$.
\end{remark}

\subsection{Classification of quasisymmetric categories}

Let $\overline{GT}=\overline{GT}(\mathbb C)$ be the
Grothendieck-Teichm\"uller semigroup, defined by Drinfeld
\cite{Dr2}. Recall that this semigroup consists of pairs
$(\lambda,f)$, where $\lambda\in \mathbb C$ and
$f(A,B)=e^{\tilde{f}(\log A,\log B)}$,
where $\tilde{f}$ is a formal Lie series satisfying some
properties.

Recall (see e.g. \cite{EK2}, Section 2.2) that the semigroup $\overline{GT}$
acts on the set of equivalence classes of quasisymmetric
categories. Namely, let $\cC$ be a quasisymmetric category.
The action of $g=(\lambda,f)$ on $\cC$
is given by preserving the abelian category structure and the
functor of tensor product, and transforming the associativity
isomorphism and braiding by the formulas
$$
\beta'=\beta\circ (\beta^2)^{\frac{\lambda-1}{2}},
$$
$$
\alpha'=\alpha \circ f(\beta^2_{12},\alpha^{-1}\circ
\beta^2_{23}\circ \alpha).
$$

\begin{remark}
Note that $\beta'$ in the above formula is well defined since
$\beta^2 -1$ is nilpotent. Namely, for any complex number $s$,
$(\beta^2)^s$ is by definition equal to $e^{s \log \beta^2}$, where
$$
\log \beta^2 = \log (1+(\beta^2 -1))=\sum_{m\ge 1} (-1)^{m-1} \frac{(\beta^2
-1)^m}{m}.
$$
\end{remark}

Recall also (\cite{Dr2}, Proposition 5.2)
that every Lie associator $\Phi$ gives
rise to a canonical 1-parameter subsemigroup $g_\Phi(\lambda)=(\lambda,
f_\Phi(\lambda))$.

We can now state our first main result.

\begin{theorem}\label{main1}
(i) Any quasisymmetric tensor category is equivalent, as a braided
tensor category, to $\mathcal{C}(G,u,t)$ for some $(G,u,t)$ with
nilpotent $t\in (S^2\g)^G$.

(ii) $\mathcal{C}(G,u,t)$ is equivalent to $\mathcal{C}(G',u',t')$
if and only if there exists a supergroup isomorphism $\phi:G\to
G'$ sending $u$ to $u'$ such that $(d\phi\ot d\phi)(t)=t'$.

(iii) Any braided unipotent tensor category is equivalent, as a
braided tensor category, to $\cC(G,t)$, where $G$ is a
prounipotent proalgebraic group, and $t\in (S^2\g)^G$. The
pair $(G,t)$ is determined uniquely up to an isomorphism.
\end{theorem}

\begin{proof}
Let $\cC$ be a quasisymmetric category.
Consider the 1-parameter family of quasisymmetric categories
$\cC(\lambda):=g_\Phi(\lambda)(\cC)$ (so $\cC(1)=\cC)$.
This family depends polynomially on $\lambda$: the
abelian category structure and the tensor product functor
do not change with $\lambda$, while the associativity isomorphism
and the braiding depend polynomially on $\lambda$.

The main point is that the category $\cC(0)$ is symmetric.
Therefore, by Deligne's theorem \cite{De2}, it is equivalent to
${\rm Rep}(G,u)$ for some $(G,u)$. Thus we can identify $\cC$ with
$\Rep(G,u)$ as an abelian category with the tensor product
functor. Then the category $\cC(\lambda)$ can be described as
follows. The braiding in this category is $P\circ e^{\lambda
t/2}$, where $P$ is the symmetry morphism of $\Rep(G,u)$ and
$t=\log(\beta^2)$, and the associativity isomorphism is
$\alpha=\Phi(\lambda t_{12},\lambda t_{23})$. Expanding the
hexagon relations in powers of $\lambda$ and taking the linear
part, we find that $t_{12,3}=t_{13}+t_{23}$, which implies that
$t\in (S^2\g)^G$. Moreover, because of the quasisymmetry
condition, $t$ is nilpotent. Thus, setting $\lambda=1$, we get
$\cC=\cC(G,u,t,\Phi)=\cC(G,u,t)$. This proves part (i) of the
theorem.

To prove part (ii), assume that we have an equivalence
of braided tensor categories $F: \cC(G,u,t,\Phi)\to \cC(G',u',t',\Phi)$.
Applying the semigroup $g_\Phi(\lambda)$ to this equivalence, we
get an equivalence $F_\lambda: \cC(G,u,\lambda t,\Phi)\to
\cC(G',u',\lambda t',\Phi)$. By setting $\lambda=0$ and using
Deligne's theorem, we can assume that $(G,u)=(G',u')$ and $F$ is
the identity functor. Then we get $t=t'$, as desired.
This proves (ii).

Part (iii) follows from (i) and (ii).

This completes the proof of the theorem.
\end{proof}

\section{Fiber functors on quasisymmetric categories}

\subsection{Unipotent fiber functors}

Let $\mathcal{C}$ be a quasisymmetric category (so
$\cC=\mathcal{C}(G,u,t)$, where $t$ is nilpotent).

\begin{definition}
Let us say that a fiber functor $F: \cC\to \Vect$ is unipotent
if it coincides with the standard one on the subcategory
$\cC_s=\Rep(G_{\rm red},u)$.
\end{definition}

Let us give a construction of unipotent fiber functors.
(In the case $t=0$, this construction appears in \cite{EG1},
Theorem 5.5.) Let $r\in \mathfrak{g}\ot \mathfrak{g}$ be a nilpotent solution
to the classical Yang-Baxter equation
$$
[r_{12},r_{13}]+[r_{12},r_{23}]+[r_{13},r_{23}]=0,
$$
such that $t=r+r_{12}$. Then by the results of \cite{EK1},
\cite{EK2}, \cite{E}, there exists a universal formula $J=J(\hbar
r)$ defining a pseudotwist killing the associator $\Phi(\hbar
t_{12},\hbar t_{23})$. Since $r$ is nilpotent, similarly to Remark
\ref{nilpo}, the series $J(\hbar r)$, when evaluated in any
product $X\otimes Y$, is in fact a polynomial in
$\hbar$. Thus, it can be evaluated at $\hbar=1$. This gives rise
to a unipotent fiber functor $F_r: \cC\to \Vect $, which is the
usual forgetful functor on $\Rep(G,u)$ with tensor structure
defined by $J(r)$.

\subsection{Classification of unipotent fiber functors}

Our second main result is the following one.

\begin{theorem}\label{main2}
If $F$ is a unipotent fiber functor on $\mathcal{C}$ then $F=F_r$
for some nilpotent $r$, and
$r$ is uniquely determined up to conjugation.
\end{theorem}

\begin{proof}
As we have already mentioned above, according to \cite{EK1,EK2,E},
the quantization $(U_\hbar(\g),R)$ of a quasitriangular Lie
bialgebra $(\g,r)$ can be obtained by twisting the enveloping
algebra $U(\g)[[\hbar]]$ by a pseudotwist $J(\hbar r)$, given by a
certain universal formula. The quantum $R$-matrix of $U_\hbar(\g)$
is then given by the universal formula
$$
R(\hbar r)=J_{21}(\hbar r)e^{\hbar(r+r_{21})/2}J(\hbar r)^{-1}=
1+\hbar r+O(\hbar^2).
$$
Thus
$$
\log R=H(\hbar r)=\hbar r+O(\hbar^2),
$$
where $H$ is some infinite series.
Since $H$ is the identity modulo higher terms, this formula can be inverted:
\begin{equation}\label{Hinv}
\hbar r=H^{-1}(\log R).
\end{equation}

Now suppose that $F$ is a unipotent fiber functor on $\cC$. Then
$B={\rm End}(F)$ is a (topological) quasitriangular Hopf algebra.
Let us now apply formula (\ref{Hinv}) to the R-matrix of $B$,
$R=P\circ F(\beta)$. Because of the unipotency of $F$, this
R-matrix is unipotent, and hence $H^{-1}(\log R)$ makes sense,
even though $H^{-1}$ is an infinite series, and there is no formal
parameter $\hbar$. So we can set $r=H^{-1}(\log R)$. Let us twist
$B$ by the twist $J(r)^{-1}$, and denote the corresponding Hopf
algebra by $B_0$. It follows from the fact that the quantization
of quasitriangular Lie bialgebras gives rise to a prop isomorphism
(Section 5 of \cite{EE}) that $B_0$ is cocommutative, and the
category $\Rep(B_0)$ with trivial symmetric structure is
equivalent to $\Rep(G,u)$ as a symmetric tensor category.
Moreover, $r\in \g\otimes \g$, and $B_0$ is equipped with a
quasitriangular co-Poisson structure defined by $r$. It is now
clear that $F=F_r$. This proves the existence part of the theorem.
The uniqueness of $r$ up to conjugation follows from the
canonicity of the above construction (see also Theorem 5.3 of
\cite{EE}).
\end{proof}

\begin{remark} Note that if $\cC$ is a unipotent category, then
the conditions of unipotency of $F$ and nilpotency of $r$
are vacuous and can be dropped.
\end{remark}

\begin{corollary}\label{tria}
Let $\cC=\Rep(G,t)$, where $G$ is a prounipotent
proalgebraic group. Then the assignment $r\to F_r$ defines a
bijection between isomorphism classes of fiber functors on $\cC$
and elements $r\in \g\otimes \g$ which satisfy the classical
Yang-Baxter equation and the condition $r+r^{21}=t$.
\end{corollary}

Note that if $G$ is a unipotent algebraic group, and $r$ is
as in Corollary \ref{tria} with $t=0$ (i.e., $r\in \wedge^2\g$)
then by a well known theorem of Drinfeld (see \cite{ES}), the
image of $r$ (regarded as a map $\g^*\to \g$) is a Lie subalgebra
$\h\subset \g$ defining a closed subgroup $H\subset G$, and
$\omega=r^{-1}$ is a nondegenerate 2-cocycle (i.e., a
left-invariant symplectic form) on $H$. Conversely, for any closed
subgroup $H\subset G$ and a left-invariant symplectic form
$\omega$ on $H$, the element $r=\omega^{-1}$, regarded as an
element of $\wedge^2\g$, is a solution of the classical
Yang-Baxter equation. Thus, Corollary \ref{tria} implies

\begin{corollary}\label{tria1}
Let $\cC=\Rep(G)$, where $G$ is a unipotent algebraic
group. Then equivalence classes of fiber functors on $\cC$ are in
bijection with conjugacy classes of pairs $(H,\omega)$, where $H$
is a closed subgroup of $G$, and $\omega$ is a left invariant
symplectic form on $H$.
\end{corollary}

This corollary is an analog, for unipotent algebraic
groups, of the classification of fiber functors on $\Rep(G)$ for
finite group $G$, due to Movshev \cite{Mo}.

\begin{remark}
We note that a classification of fiber functors
on $\Rep(G,u)$ for a general
affine algebraic supergroup $G$ is unknown.
It is clear that to obtain such a classification,
it would be sufficient to do so for $G=GL(m|n)$, which is equivalent
to classifying unitary solutions of the quantum Yang-Baxter
equation which are invertible and skew-invertible. This problem
is open even for $n=0$ (starting from $m=4$).
\end{remark}

\subsection{Classification of coconnected coquasitriangular Hopf
algebras}

\begin{definition} (see e.g. \cite{S})
A Hopf algebra $A$ is called {\it coconnected} if
every simple comodule over $A$ is trivial.
\end{definition}

\begin{corollary}\label{tria2}
Let $A$ be a coconnected coquasitriangular Hopf algebra. Then $A$
is obtained by twisting the product of the function algebra
${\mathcal O}(G)$ of a prounipotent proalgebraic group $G$
by a pseudotwist $J=J(r)$, where $r\in \g^{\otimes 2}$ is a
solution of the classical Yang-Baxter equation, such that
$t=r+r^{21}$ is $G$-invariant. Moreover, the pair $(G,r)$ is
determined by $A$ up to an isomorphism.
\end{corollary}

\begin{proof}
Follows from
Theorem \ref{main1}(iii) and Corollary \ref{tria},
since the comodule category of a
coconnected coquasitriangular Hopf algebra is a unipotent braided
tensor category.
\end{proof}

\section{Classification of coconnected Hopf algebras}

In this section we will give a classification of
coconnected Hopf algebras.

\subsection{Construction of coconnected Hopf algebras}

Let $G$ be a prounipotent Poisson proalgebraic group, and
$\g$ its Lie algebra. Then $\g$ has a Lie bialgebra structure
$\delta$ which determines the Poisson-Lie structure on $G$.
Consider the lower central series filtration $\mathfrak{g}_\bullet$
on $\mathfrak{g}$: $\mathfrak{g}_0=\mathfrak{g}$ and
$\mathfrak{g}_{i+1}=[\mathfrak{g}_i,\mathfrak{g}]$. Then the
bracket $[,]$ of $\g$ has degree $1$ with respect to this
filtration, while the cobracket $\delta$ has degree $0$ (since
$\delta([a,b])=[a\ot 1+1\ot a,\delta(b)]+ [\delta(a),b\ot 1+1\ot
b]$).

Consider now the Etingof-Kazhdan quantization ${\mathcal
O}_\hbar(G)$ of the Poisson group $G$, \cite{EK3}. This is the
space ${\mathcal O}(G)[[\hbar]]$ with certain product and
coproduct deforming the standard Hopf algebra structure on
${\mathcal O}(G)$. These product and coproduct are given by some
universal formulas (infinite series) in terms of the bracket and
cobracket of $\g$. Because the degree of the bracket is 1 and of the
cobracket is zero, these infinite series terminate, and thus the
formal quantization of $G$ is actually defined over the polynomials
$\mathbb C[\hbar]$, which means that the Hopf algebra ${\mathcal
O}_\hbar(G)$ has a lattice ${\mathcal O}_\hbar(G)_{\rm pol}$ over
$\mathbb C[\hbar]$. Specializing this lattice to $\hbar=1$, we get
a Hopf algebra $A=A(G,\delta)$ over $\mathbb C$. As a vector
space, it coincides with ${\mathcal O}(G)$.

\begin{lemma}\label{coco}
The Hopf algebra $A$ is coconnected.
\end{lemma}

\begin{proof}
Let $A_0:={\mathcal O}(G)$, and let $C_\bullet(A_0)$ be the
coradical filtration of $A_0$. Since $A_0$ is coconnected,
$C_0(A_0)=\mathbb C$. It is easy to see that the product in $A_0$
preserves this filtration, while the Poisson bracket preserves it
in the strict sense, i.e., decreases the filtration degree by $1$.
This implies that the product and coproduct in $A$ preserve the
coradical filtration (as they are obtained from the product,
coproduct, and Poisson bracket of $A_0$ by a universal formula).
In particular, the coradical filtration $C_\bullet(A)$ of $A$
coincides with $C_\bullet(A_0)$. This implies that $A$ is
coconnected.
\end{proof}

\subsection{Classification of coconnected Hopf algebras}

Our third main result is the following theorem.

\begin{theorem}\label{main3}
Any coconnected Hopf algebra over $\mathbb C$ is of the form
$A(G,\delta)$, and $(G,\delta)$ is determined up to an
isomorphism.
\end{theorem}

\begin{proof}
The results of \cite{EK2},\cite{EE} imply that the formulas
expressing the product $\mu$ and coproduct
$\Delta$ of ${\mathcal O}_\hbar(G)$ in terms of the product $\mu_0$, coproduct
$\Delta_0$, and bracket $\delta$ of ${\mathcal O}(G)$ are
invertible. Let $A$ be a coconnected Hopf algebra,
and let us apply the inverse formulas
at $\hbar=1$ to introduce a commutative
product $\mu_0$, coproduct $\Delta_0$, and bracket $\delta$ on
$A$. As before, the formulas make sense for $\hbar=1$
because series terminate due to the pronilpotency of $\g$.
In this way we get a commutative Poisson-Hopf algebra
$A_0$, which coincides with $A$ as a vector space.

Let $C_\bullet(A)$ be the coradical filtration of $A$. Since $A$
is coconnected, $C_0(A)=\mathbb C$. Also, this filtration is fixed
by the product and coproduct $\mu,\Delta$ in $A$. Hence it is
fixed by the new coproduct $\Delta_0$, as $\Delta_0$ expresses via
$\mu,\Delta$ by a universal formula. Thus $C_\bullet(A)$ coincides
with the coradical filtration $C_\bullet(A_0)$ of $A_0$. Hence
$A_0$ is coconnected, i.e., $A_0={\mathcal O}(G)$. It is now clear
that $A=A(G,\delta)$. The uniqueness of $G,\delta$ is clear from
the canonicity of this construction. The theorem is proved.
\end{proof}

\begin{remark}
A classical theorem of Kostant states that
a coconnected {\it cocommutative} Hopf
algebra over $\mathbb C$ is the same thing as an enveloping
algebra of a Lie algebra.
This fits with Theorem \ref{main3} as follows. It is easy
to show that the cocommutativity condition of $A(G,\delta)$ is
equivalent to the condition that $G$ is abelian. Thus $G=\g$ is a
provector space. Denote the topological dual space $\g^*$ by $L$;
it is an ordinary (discrete) vector space, possibly infinite dimensional.
Then $\delta$ gives rise to a Lie bracket on $L$, and it is easy
to see that $A(G,\delta)=U(L)$.
\end{remark}

\begin{remark}
In categorical terms, Theorem \ref{main3} provides a
classification of unipotent tensor categories with a fiber
functor; namely, they are categories of finite dimensional
comodules over Hopf algebras $A(G,\delta)$.
\end{remark}

\begin{remark}
If $\cC$ is a {\it finite} unipotent tensor category, then
it is shown in \cite{EO}, Section 2.10, that
$\cC$ is the category of vector spaces. This is a special case of
Theorem \ref{main3}, since in the finite case we must have $G=1$.
\end{remark}

\begin{remark}
One may hope that along these lines one should be able to obtain
a classification of general unipotent tensor categories; namely,
one could expect that their equivalence classes are in a natural
bijection with equivalence classes of prounipotent groups
$G$ with a Lie quasibialgebra structure on ${\rm Lie}(G)$.
Unfortunately, this remains out of reach, as a quantization
theory of Lie quasibialgebras is still unavailable.
\end{remark}

\begin{remark}
We have obtained several theorems giving a classification of
various kinds of tensor categories in terms of affine proalgebraic
supergroups with some additional data. We note that in these
theorems, an additional requirement that the relevant affine
proalgebraic supergroup be actually {\it algebraic} is equivalent to the
requirement that the corresponding category $\cC$ be {\it finitely
$\otimes$-generated} (\cite{De2}), i.e., there exists an object $X$
in $\mathcal{C}$ such that every object in $\mathcal{C}$ is a
subquotient of a direct sum of tensor powers of $X$.
\end{remark}

\end{document}